\documentclass[a4paper,12pt]{article}
\usepackage{comment}
\usepackage{cite}
\usepackage{amsmath}
\usepackage{amssymb}
\usepackage{amsfonts}
\usepackage[T1]{fontenc}
\usepackage[utf8]{inputenc}
\usepackage{graphicx}
\usepackage{fancyhdr}
\usepackage{float}
\usepackage{xcolor}
\usepackage{authblk}
\usepackage{mathrsfs}
\usepackage{empheq}
\usepackage[hyphens]{url}
\usepackage{hyperref} 
\usepackage[]{breakurl}

\usepackage{geometry}
\begin{document}

\title{An exact series expansion for the Dottie number}

\author[$\dagger$]{Jean-Christophe {\sc Pain}$^{1,2,}$\footnote{jean-christophe.pain@cea.fr}\\
\small
$^1$CEA, DAM, DIF, F-91297 Arpajon, France\\
$^2$Universit\'e Paris-Saclay, CEA, Laboratoire Mati\`ere en Conditions Extr\^emes,\\ 
91680 Bruy\`eres-le-Ch\^atel, France
}

\maketitle

\begin{abstract}
In this article, an exact series expansion for the Dottie number (solution of the equation $\cos x=x$) is presented. Its derivation consists in combining the Kaplan representation of the Dottie number as a series in odd powers of $\pi$, with a series for $\pi^{2k+1}$ obtained using the Fa\`a di Bruno formula. 
\end{abstract}

\section{Introduction}

The Dottie number is the name given by Kaplan \cite{Kaplan2007} to the unique real root of the equation $\cos x=x$ (namely, the unique real fixed-point of the cosine function) \cite{Fettis1976,Miller1980,Azarian2007,Cottanceau2008,Salov}. Actually, the name ``Dottie'' would be, according to the legend, referring to a particular French professor who noticed that whenever she typed a number into her calculator and hit the cosine button repeatedly, the result always converged to this value. The number appears in numerous elementary works on algebra already \cite{Bertrand1865}. It is also known sometimes simply as the cosine constant, cosine superposition constant, iterated cosine constant, or cosine fixed point constant. Arakelian has used the Armenian small letter ``ayb'' (first letter in the Armenian alphabet) to denote this constant \cite{Arakelian1981,Arakelian1995}.

This root is a simple nontrivial example of an universal attracting fixed point. It is also transcendental as a consequence of the Lindemann-Weierstrass theorem (see Appendix A) and finds applications in mathematics, for instance in the Bertrand semi-circle problem (see Appendix B), or in celestial physics, in connection with the Kepler equation (see Appendix C). Let us denote $\mathscr{D}$ the Dottie number. Its first 32 significant digits are

\begin{equation}
\mathscr{D}=0.739085 133215 160641 655312 087673 87...    
\end{equation}
and one has
\begin{equation}
\sin(\mathscr{D}) \approx 0.673612029183...
\end{equation}
as well as
\begin{equation}
\tan(\mathscr{D}) \approx 0.911413312094....
\end{equation}
Simpe approximations of the Dottie number are provided in Appendix D. The following explicit form was obtained \cite{Gaidash2022}:
\begin{equation}
\mathscr{D}=\sqrt{1-\left(2I_{1/2}^{-1}\left(\frac{1}{2},\frac{3}{2}\right)-1\right)^2},
\end{equation}
where $I_{1/2}^{-1}$ denotes the inverse of the regularized incomplete Beta function
\begin{equation}
I_x(a,b)=\frac{B(x;a,b)}{\mathrm{B}(a,b)},
\end{equation}
where
\begin{equation}
B(x;a,b)=\int_0^xt^{a-1}(1-t)^{b-1}dt
\end{equation}
is the incomplete Beta function and $\mathrm{B}(a,b)=B(1;a,b)$ the usual Beta function. The Dottie number can be expressed as
\begin{equation}\label{eq7}
\mathscr{D}=\sum_{n=0}^{\infty}b_n\pi^{2n+1},
\end{equation}
and the latter formula will be the starting point of the present work.

It is worth mentioning that the Dottie number can also be expressed as an Engel series \cite{EngelOEIS,Broukhis,Engel1913,Erdos1958}.

\section{Explicit form of the coefficients $a_n$ in the Kaplan series}

Let us set $f(x)=x-\cos x$, and $g(x)=f^{-1}(x)$ its inverse (reciprocal) function, which has \emph{a priori} no explicit form. The Taylor-series expansion of $g$ reads
\begin{equation}
g(x)=\sum_{n=0}^{\infty}g^{(n)}(x_0)\frac{(x-x_0)^n}{n!}.
\end{equation}
Since $f(\pi/2)=\pi/2$ and $g(\pi/2)=\pi/2$, it is convenient to choose $x_0=\pi/2$ and we get
\begin{equation}
g(x)=\sum_{n=0}^{\infty}g^{(n)}\left(\frac{\pi}{2}\right)\frac{1}{n!}\left(x-\frac{\pi}{2}\right)^n.
\end{equation}
The Dottie number is therefore
\begin{equation}
\mathscr{D}=g(0)=\sum_{n=0}^{\infty}g^{(n)}\left(\frac{\pi}{2}\right)\frac{\left(-\pi\right)^n}{2^nn!}.
\end{equation}
The multiple derivatives $g^{(n)}$ of $g$ can be deduced from the multiple derivatives $f^{(n)}$ of $f$ using the Fa\`a di Bruno formula. One has in particular $f'(x)=1+\sin x$, $f''(x)=\cos x$ and $\forall~ n>1$:
\begin{equation}
f^{(n)}(x)=\frac{d^{n-2}}{dx^{n-2}}\cos x
\end{equation}
yielding the values of table \ref{tab1}.
\begin{table}
\begin{center}
\begin{tabular}{ccccccccccc}\hline
$n$ & 0 & 1 & 2 & 3 & 4 & 5 & 6 & 7 & 8 & 9\\\hline
$f^{(n)}\left(\frac{\pi}{2}\right)$ & $\frac{\pi}{2}$ & 2 & 0 & -1 & 0 & 1 & 0 & -1 & 0 & 1\\ \hline
\end{tabular}
\end{center}
\caption{Values of multiple derivatives (orders from zero to nine) of function $f:x\longmapsto\cos x$, for $x=\pi/2$.}\label{tab1}
\end{table}

\begin{table}
\begin{center}
\begin{tabular}{ccccccccccc}\hline
$n$ & 0 & 1 & 2 & 3 & 4 & 5 & 6 & 7 & 8 & 9\\\hline
$g^{(n)}\left(\frac{\pi}{2}\right)$ & $\frac{\pi}{2}$ & $-\frac{1}{768}$ & 0 & $-\frac{1}{61440}$ & 0 & $-\frac{43}{165150720}$ & 0 & $-\frac{233}{47563407360}$ & 0 & $-\frac{60623}{669692775628800}$\\ \hline
\end{tabular}
\end{center}
\caption{Values of multiple derivatives (orders from zero to nine) of function $g: x\longmapsto f^{-1}(x)$, where $f^{-1}$ is the reciprocal function of $f$ ($f^{-1}(f(x))=f(f^{-1}(x))=x$) for $x=\pi/2$.}\label{tab2}
\end{table}
The first terms of the series are thus \cite{Ozaner}:
\begin{eqnarray}
\mathscr{D}&=&\pi/4-\pi^3/768-\pi^5/61440-43\pi^7/165150720-233\pi^9/47563407360\nonumber\\
& &-\frac{60623}{669692775628800}\pi^{11}.
\end{eqnarray}
The above derivation of the $a_n$ coefficients can be reformulated invoking Lagrange inversion theorem \cite{Apostol2000,Johnson2002}. If $y=f(x)$, where $f(0)=0$ and $f'(0)\ne  0$, then
\begin{equation}
x=\sum_{n=1}^{\infty}\frac{y^n}{n!}\left[\frac{d^{n-1}}{dx^{n-1}}\left(\frac{x}{f(x)}\right)^n\right]_{x=0}.
\end{equation}
In fact, one has 
\begin{equation}\label{an}
\mathscr{D}=\frac{\pi}{2}+\sum_{n~\mathrm{odd}}a_n\pi^n,
\end{equation}
where
\begin{equation}
a_n=\frac{1}{n!2^n}\lim_{x\rightarrow\frac{\pi}{2}}\frac{\partial^{n-1}}{\partial x^{n-1}}\left[\frac{\cos x}{x-\frac{\pi}{2}}-1\right]^{-n}.
\end{equation}
Comparing with Eq. (\ref{eq7}), we have $b_1=a_1+\pi/2$, $b_n=a_n$ for $n\geq 3$ and $n$ odd, and $b_n=0$ if $n$ even. The idea is to use the Fa\`a di Bruno formula to calculate 
\begin{equation}
\left(\frac{1}{h(x)}\right)^{(k)},
\end{equation}
with
\begin{equation}
h(x)=\left[\frac{\cos x}{x-\frac{\pi}{2}}-1\right]^{n}.
\end{equation}
We have
\begin{equation}
\left(\frac{1}{h(x)}\right)^{(k)}=\frac{k!}{h^{k+1}(x)}\sum_{\{m_i\}}\frac{(-1)^{k-m_0}(k-m_0)!}{\prod_{i=1}^k(i!)^{m_i}m_i!}\prod_{i=0}^n\left(h^{(i)}(x)\right)^{m_i},
\end{equation}
with $m_0=m_2+2m_3+\cdots+(n-1)m_k$ and $n=m_0+m_1+m_2+\cdots+m_n=m_1+2m_2+\cdots+km_k$.
We have
\begin{equation}
h^{(i)}(x)=\sum_{\{p_l\}/p_1+p_2+\cdots+p_i=i}\binom{i}{p_1,\cdots,p_i}\prod_{j=1}^i\frac{\partial^{p_j}}{\partial x^{p_j}}\left(\frac{\cos x}{x-\frac{\pi}{2}}-1\right).
\end{equation}
Using
\begin{equation}
\frac{\partial^{p_j}}{\partial x^{p_j}}\left(\frac{\cos x}{x-\frac{\pi}{2}}\right)=\sum_{u=0}^{p_j}\binom{p_j}{u}\cos^{(u)}(x)\times\frac{(-1)^{p_j-u}}{\left(x-\frac{\pi}{2}\right)^{p_j-u+1}}(p_j-u)!
\end{equation}
and since
\begin{equation}
\cos^{(p)}(x)=\cos\left(x+p\frac{\pi}{2}\right),
\end{equation}
we get, with
\begin{equation}
\left(\frac{1}{h(x)}\right)^{(n-1)}=\frac{(n-1)!}{h^{n}(x)}\sum_{\{m_i\}}\frac{(-1)^{n-1-m_0}(n-1-m_0)!}{\prod_{i=1}^{n-1}(i!)^{m_i}m_i!}\prod_{i=0}^{n-1}\left(h^{(i)}(x)\right)^{m_i},
\end{equation}
the final formula
\begin{empheq}[box=\fbox]{align}\label{jc}
a_n&=\frac{1}{n!2^n}\lim_{x\rightarrow\frac{\pi}{2}}\frac{(n-1)!}{\left(\frac{\cos x}{x-\frac{\pi}{2}}-1\right)^{n^2}}\sum_{\{m_i\}}\frac{(-1)^{n-1-m_0}(n-1-m_0)!}{\prod_{i=1}^{n-1}(i!)^{m_i}m_i!}\nonumber\\
& \times\prod_{i=0}^{n-1}\left[\sum_{\{p_l\}/p_1+p_2+\cdots+p_i=i}\binom{i}{p_1,\cdots,p_i}\prod_{j=1}^i\left\{\sum_{u=0}^{p_j}\frac{p_j!}{u!}\cos\left(x+u\frac{\pi}{2}\right)\frac{(-1)^{p_i-u}}{(x-\frac{\pi}{2})^{p_i-u+1}}\right\}\right]^{m_i}.\nonumber\\
& 
\end{empheq}

\section{Series expansion of $\pi^{n}$, with $n$ an odd integer}

We recently published a closed-form expression for $\pi^{2+k}$ \cite{Pain2022}:
\begin{equation}\label{newpik}
\pi^{2+k}\sum_{\substack{\sum_{j=1}^kjm_j=k}}\frac{\left(\sum_{j=1}^km_j\right)!~2^k\prod_{j=1}^k\left[\cos\left(2\pi x+j\frac{\pi}{2}\right)\right]^ {m_j}}{2^{\sum_{j=1}^km_j}\left(\prod_{j=1}^km_j!\right)\left[\sin^2(\pi x)\right]^{\sum_{j=1}^km_j+1}\prod_{j=1}^k(j!)^{m_j}}=\sum_{n=-\infty}^{\infty}\frac{(-1)^k(k+1)}{(x-n)^{k+2}},
\end{equation}
which can be put in the form
\begin{equation}\label{pi}
\pi^{2+k}=\frac{(-1)^k(k+1)}{2^k\mathscr{A}_k(x)}\sum_{n=-\infty}^{\infty}\frac{1}{(x-n)^{k+2}},
\end{equation}
where we have set
\begin{equation}
\mathscr{A}_k(x)=\sum_{\left\{m_i\right\}}\frac{\mathscr{S}_k!}{2^{\mathscr{S}_k}\left(\prod_{j=1}^km_j!\right)\left[\sin^2(\pi x)\right]^ {\mathscr{S}_k+1}}\times\prod_{j=1}^k\frac{\left[\cos\left(2\pi x+j\frac{\pi}{2}\right)\right]^{m_j}}{(j!)^{m_j}},
\end{equation}
and
\begin{equation}
\mathscr{S}_k=\sum_{j=1}^km_j
\end{equation}
with
\begin{equation}
\sum_{j=1}^kjm_j=k.
\end{equation}
Formula (\ref{newpik}) was obtained starting from the relation \cite{Borwein1987}:
\begin{equation}
\pi~\mathrm{cotan}(\pi x)-\pi~\mathrm{cotan}(\pi a)=\sum_{n=-\infty}^{\infty}\left(\frac{1}{x-n}-\frac{1}{a-n}\right)=\sum_{n=-\infty}^{\infty}\frac{(a-x)}{(x-n)(a-n)},
\end{equation}
and gives for instance
\begin{equation}\label{euler}
\pi^3\left[\cot(\pi x)~\mathrm{cosec}^2(\pi x)\right]=\sum_{n=-\infty}^{\infty}\frac{1}{(x-n)^3}
\end{equation}
and
\begin{equation}
\pi^5\cot(\pi x)\mathrm{cosec}^2(\pi x)\left[\mathrm{cosec}^2(\pi x)-\frac{1}{3}\right]=\sum_{n=-\infty}^{\infty}\frac{1}{(x-n)^5}.
\end{equation}
\section{Series expansion of the Dottie number}

Combining Eq. (\ref{jc}) with Eqs. (\ref{an}) and (\ref{pi}) yields the following expression for the Dottie number:
\begin{empheq}[box=\fbox]{align}\label{fina}
\mathscr{D}&=\frac{\pi}{2}+\sum_{n~\mathrm{odd}\geq 3}a_n\pi^n\frac{1}{n!2^n}\lim_{x\rightarrow\frac{\pi}{2}}\frac{(n-1)!}{\left(\frac{\cos x}{x-\frac{\pi}{2}}-1\right)^{n^2}}\sum_{\{m_i\}}\frac{(-1)^{n-1-m_0}(n-1-m_0)!}{\prod_{i=1}^{n-1}(i!)^{m_i}m_i!}\nonumber\\
& \times\prod_{i=0}^{n-1}\left[\sum_{\{p_l\}/p_1+p_2+\cdots+p_i=i}\binom{i}{p_1,\cdots,p_i}\prod_{j=1}^i\left\{\sum_{u=0}^{p_j}\frac{p_j!}{u!}\cos\left(x+u\frac{\pi}{2}\right)\frac{(-1)^{p_i-u}}{\left(x-\frac{\pi}{2}\right)^{p_i-u+1}}\right\}\right]^{m_i}\nonumber\\
& \times\frac{(-1)^{n}(n-1)}{2^{n-2}\mathscr{A}_{n-2}(x)}\sum_{p=-\infty}^{\infty}\frac{1}{(x-p)^{n}},
\end{empheq}
with
\begin{equation}
\mathscr{A}_{n-2}(x)=\sum_{\left\{m_i\right\}}\frac{\mathscr{S}_{n-2}!}{2^{\mathscr{S}_{n-2}}\left(\prod_{j=1}^{n-2}m_j!\right)\left[\sin^2(\pi x)\right]^{\mathscr{S}_{n-2}+1}}\times\prod_{j=1}^{n-2}\frac{\left[\cos\left(2\pi x+j\frac{\pi}{2}\right)\right]^{m_j}}{(j!)^{m_j}},
\end{equation}
and
\begin{equation}
\mathscr{S}_{n-2}=\sum_{j=1}^{n-2}m_j
\end{equation}
together with
\begin{equation}
\sum_{j=1}^{n-2}jm_j=n-2.
\end{equation}
We did not manage to simplify expression (\ref{fina}) further.

\section{Conclusion}

In this article, an exact series expansion for the Dottie number (solution of the equation $\cos x=x$) was presented. It relies on the combination of the Kaplan series, which expresses the Dottie number as a series in odd powers of $\pi$, with a series for $\pi^{2k+1}$ ($k\geq 1$) recently obtained using the Fa\`a di Bruno formula. The resulting formula, however, is rather cumbersome and we hope that the present work will stimulate the derivation of additional explicit formulas or closed-forms of the Dottie number.

\appendix

\section*{Appendix A: On the transcendance of the Dottie number}\label{appA}

An algebraic number $x$ is a number satisfying some polynomial equation $a_nx^n + a_{n-1}x^{n-1} + \cdots + a_2x^2 + a_1x + a_0 = 0$, where each $a_i$ is rational. This definition is valid for complex numbers as well. It is known that if $x$ and $y$ are algebraic, it is also the case for $x+y$, $x-y$, and $xy$ \cite{Niven}. Furthermore, it can easily be shown that if $x$ is algebraic then so is the square root of $x$ (since $x^n = (\sqrt{x})^{2n}$). As a consequence of the Lindemann–Weierstrass theorem \cite{Lindemann1882a,Lindemann1882b,Weierstrass1885,Niven}, if $x \ne 0$ is an algebraic number then $e^x$ is transcendental. Let us assume $\mathscr{D}$ is algebraic. We have
\begin{equation}
e^{i\mathscr{D}}=\cos(\mathscr{D})+i\sin(\mathscr{D}).
\end{equation}
From $\sin^2x + \cos^2x = 1$ and $\mathscr{D} = \cos \mathscr{D}$, the latter equation simplifies to
\begin{equation}\label{exp}
e^{i\mathscr{D}}=\mathscr{D}+i\sqrt{1-\mathscr{D}^2}.
\end{equation}
Since the right side of Eq. (\ref{exp}) involves addition, subtraction, multiplication, and a square root, the expression is algebraic. Therefore if $\mathscr{D}$ is algebraic, so must be $e^{i\mathscr{D}}$.

Since $i$ is algebraic and $\mathscr{D}$ is assumed algebraic, $i\mathscr{D}$ is algebraic too and by the Lindemann-Weierstrass theorem, $e^{i\mathscr{D}}$ is transcendental. The contradiction with the result above means that $i\mathscr{D}$ can not be algebraic. In other words, the Dottie number is transcendental. However, it might still be represented in terms of $\pi$, oas we have seen above.

\section*{Appendix B: Connection with the Bertrand semi-circle}\label{appB}

\begin{figure}
\begin{center}
\includegraphics[width=15cm]{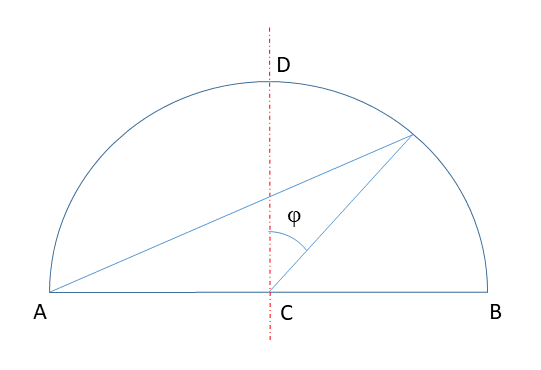}
\end{center}
\caption{Bertrand's semi-circle: the angle $\phi$ and the points A, B, C and D.}\label{fig1}
\end{figure}

\begin{figure}
\begin{center}
\includegraphics[width=15cm]{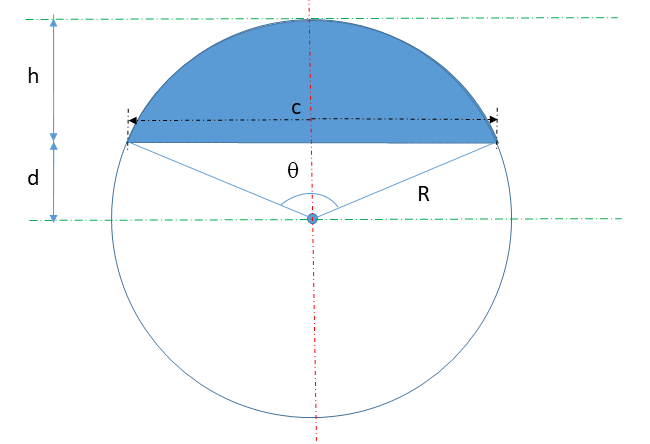}
\end{center}
\caption{Representation of Bertrand's (semi-)circle for the calculation of the area above the chord.}\label{fig2}
\end{figure}

Bertrand semi-circle problem consists in finding the point $M$ (see figure \ref{fig1}), \emph{i.e.} the angle $\phi$, such as the segment $[AM]$ splits the half-disk in two parts of equal area. Let us first calculate the area of the blue region (chord) of the disk represented in figure \ref{fig2}.

The area above the chord $[AM]$ reads
\begin{equation}
A=\pi R^2\frac{\theta}{2\pi}-A_1,
\end{equation}
where $A_1$ is the area of the triangle
\begin{equation}
A_1=\frac{c}{2}\times R\cos\left(\frac{\theta}{2}\right)
\end{equation}
with
\begin{equation}
c=2R\sin\left(\frac{\theta}{2}\right).
\end{equation}
The area of the triangle is therefore
\begin{equation}
A_1=\frac{R^2}{2}\sin\theta,
\end{equation}
and thus
\begin{equation}
A=\pi R^2\frac{\theta}{2\pi}-\frac{R^2}{2}\sin\theta=\frac{R^2}{2}\left(\theta-\sin\theta\right).
\end{equation}
The triangle is split into two parts of equal area if, setting $\theta=\pi/2+\phi$:
\begin{equation}
\frac{R^2}{2}\left[\frac{\pi}{2}+\phi-\sin\left(\frac{\pi}{2}+\phi\right)\right]=\frac{\pi R^2}{2}-\frac{R^2}{2}\left[\frac{\pi}{2}+\phi-\sin\left(\frac{\pi}{2}+\phi\right)\right]
\end{equation}
yielding
\begin{equation}
\phi=\cos\phi.
\end{equation}

\section*{Appendix C: Connection with the Kepler equation}\label{appC}

In astronomy, Kepler's equation, which relates various geometric properties of the orbit of a body subject to a central force, reads
\begin{equation}
M=E-e\sin E,
\end{equation}
$M$ being the mean anomaly, $E$ the eccentric anomaly and $e$ the eccentricity of the ellipse \cite{Rax2020}. Setting $x=a(1-e)$, $y=0$ at $t=t_0$, for $M=n(t-t_0)$, Kepler equation provides the corresponding parametrization of the ellipse
\begin{equation}
x=a\left(\cos E-e\right),\;\;\;\; y=b\sin E.
\end{equation}
The Dottie Number satisfies the ``equal area swept out in equal time'' condition (the second Kepler law) at the quarter period for $e = 1$ (parabola case). \emph{i.e.},
\begin{equation}
E-\sin(E)=\frac{\pi}{2}.
\end{equation}
Kepler's Equation reduces to
\begin{equation}
\cos(\sin(E))=\sin(E),
\end{equation}
which corresponds to the definition of the Dottie Number, where $\sin(E)=\mathscr{D}$.
The connection stems from the fact that Kepler's Equation can be expressed in terms of a Bessel function of the first kind ($J_k$):
\begin{equation}
E=M+2\sum_{n=1}^{\infty}\frac{\sin (nM)}{n}J_n(ne),
\end{equation}
yielding
\begin{equation}
\mathscr{D}=2\sum_{n=0}^{\infty}\left(\frac{J_{4n+1}(4n+1)}{4n+1}-\frac{J_{4n+1}(4n+3)}{4n+3}\right).
\end{equation}
Using the integral representation of Bessel function of the first kind

\begin{equation}
J_n(z)=\frac{1}{\pi}\int_0^{\pi}\cos(tn-z\sin t)dt,
\end{equation}
one gets \cite{kap,Integral}:
\begin{equation}
\mathscr{D}=\frac{2}{\pi}\int_0^{\pi}\left(\sum_{n=0}^{\infty}\frac{\sin[(4n+1)u]\sin[(4n+1)\sin u]}{4n+1}-\sum_{n=0}^{\infty}\frac{\sin[(4n+3)u]\sin[(4n+3)\sin u]}{4n+3}\right),
\end{equation}
which is a Kapteyn series \cite{Kapteyn1893}, \emph{i.e.} a series expansion of analytic functions on a domain in terms of the Bessel function of the first kind \cite{Watson1944}.

\section*{Appendix D: Approximations of the Dottie number $\mathscr{D}$}\label{appD}

In this Appendix, we just mention a few approximants of the Dottie number. Taking the tangent of the cosine function in $\left(\pi/4,\sqrt{2}/2\right)$ yields
\begin{equation}
\mathscr{D}\approx\frac{4+\pi}{4+4\sqrt{2}}\approx 0.7395361335...
\end{equation}
for which only the first three digits after the comma are correct. Broukhis proposed the following approximant \cite{Broukhis}:
\begin{equation}
\mathscr{D}=\left(\frac{\pi}{60}\right)^{\frac{1}{13}}\approx 0.7390853722...
\end{equation}
The first six digits after the comma are correct. This form is due to Hammond \cite{Hammond}:
\begin{equation}
\mathscr{D}\approx\Gamma^{-1}\left(e^{\pi/3}-\ln 5\right)\approx 0.7390851307...
\end{equation}
and the first eight digits after the comma are correct. 

\end{document}